\documentclass[10pt]{article}
\usepackage[tbtags]{amsmath}
\usepackage{epsfig,amstext,amssymb,amsthm,latexsym}
\usepackage{amssymb,color}

\definecolor{c20}{rgb}{0.,0.7,0.}
\definecolor{c30}{rgb}{0.,0.,1.}
\definecolor{c40}{rgb}{1,0.1,0.7}
\definecolor{c50}{rgb}{1,0,0}

\def\coEH#1{\textcolor{c50}{#1}}

\def\crA#1{\textcolor{c40}{#1}}

\def\crA#1{#1}

\catcode`\@=11 \@addtoreset{equation}{section} \catcode`\@=12

\allowdisplaybreaks[4] 
\newcommand{\nwc}{\newcommand}
\nwc{\COM}[1]{}

\nwc{\vs}[1]{\vskip #1 cm}

\newtheorem{theo}{Theorem}[section]
\newtheorem{sat}[theo]{Proposition}
\newtheorem{de}[theo]{Definition}
\newtheorem{lem}[theo]{Lemma}

\newtheorem{korr}[theo]{Corollary}
\newtheorem{remark}[theo]{Remark}
\newtheorem{exxa}[theo]{Example}

\newcommand{\nelem}[1]{{Lemma \ref{#1}}}

\newcommand{\netheo}[1]{{Theorem \ref{#1}}}

\newcommand{\kb}[1]{\boldsymbol{#1}}
\newcommand{\vk}[1]{\kb{#1}}
\newcommand{\vek}[1]{\mbox{\rm$\vk{#1}$}}
\def\FRE{\mbox{Fr\'{e}chet }}

\newcommand{\ve}{\varepsilon}

\newcommand{\abs}[1]{\lvert #1 \rvert}
\newcommand{\Abs}[1]{ \Bigl \lvert #1 \Bigr \rvert}

\newcommand{\E}[1]{\mbox{\rm$\vk{E}$}\{#1\}}

\newcommand{\pk}[1]{\mbox{\rm$\vk{P}$} \{#1\} }

\newcommand{\R}{\!I\!\!R}

\newcommand{\N}{\!I\!\!N}
\newcommand{\inr}{\in \R}

\newcommand{\inn}{\in \N}

\newcommand{\ldot}{,\ldots,}

\newcommand{\limit}[1]{\lim_{#1 \to   \infty}}

\newcommand{\todis}{\stackrel{d}{\to}}
\newcommand{\toprob}{ \stackrel{p}{\to}}

\newcommand{\equaldis}{\stackrel{d}{=}}

\newcommand{\BQN}{\begin{eqnarray}}
\newcommand{\EQN}{\end{eqnarray}}
\newcommand{\BQNY}{\begin{eqnarray*}}
\newcommand{\EQNY}{\end{eqnarray*}}
\newcommand{\BS}{\begin{sat}}
\newcommand{\ES}{\end{sat}}
\newcommand{\BL}{\begin{lem}}
\newcommand{\EL}{\end{lem}}
\newcommand{\BT}{\begin{theo}}
\newcommand{\ET}{\end{theo}}
\newcommand{\BK}{\begin{korr}}
\newcommand{\EK}{\end{korr}}
\newcommand{\BD}{\begin{de}}
\newcommand{\ED}{\end{de}}
\newcommand{\BIT}{\begin{itemize}}
\newcommand{\EIT}{\end{itemize}}
\newcommand{\BDI}{\begin{description}}
\newcommand{\EDI}{\end{description}}

\newcommand{\BEX}{\begin{exxa}}
\newcommand{\EEX}{\end{exxa}}
\def\o{\overline}

\newcommand{\QED}{\hfill $\Box$}
\newcommand{\IF}{\infty}
\def\kal#1{{\cal{ #1}}}

\newcommand{\prooftheo}[1]{ \textsc{Proof of Theorem} \ref{#1} }

\newcommand{\prooflem}[1]{\textsc{Proof of Lemma} \ref{#1}}
\newcommand{\proofkorr}[1]{\textsc{Proof of Corollary} \ref{#1}}

\def\eqdis{\equaldis}

\def\overlineALI{\alpha}
\def\overlineAL{\beta}
\def\barALI{\alpha}
\def\barAL{\beta}

\def\IDo{\coEH{p}}

\def\HAL{H_{\alpha, \lambda}}
\def\HAB{H_{\alpha, \beta}}
\def\BARHAL{\overline H_{\alpha, \lambda}}
\def\BARHAB{ \overline H_{ \alpha, \beta}}
\def\HABD{h_{\alpha,\beta}}
\def\OPJ{ \kal{J}}

\def\bab{B_{\alpha,\beta}}

\begin{document}

\centerline{\Large Tail Asymptotics under Beta Random Scaling}

        \vskip 1.8 cm
        \centerline{\large ENKELEJD HASHORVA \footnote{
Department of Mathematical  Statistics and Actuarial Science, University of Bern, Sidlerstrasse 5,
CH-3012 Bern, Switzerland, hashorva@stat.unibe.ch}
and
ANTHONY G. PAKES \footnote{
School of Mathematics and Statistics, University of Western Australia,
35 Stirling Highway, Crawley, W. A., 6009, Australia, E-mail: pakes@maths.uwa.edu.au}
}

      \centerline{\today{}}

 {\bf Abstract:} Let $X,Y,B$ be three independent
 random variables such that $X$ has the same distribution function as $Y B$.
 Assume that $B$ is a beta random variable with positive
 parameters $\alpha,\beta$ and $Y$ has distribution function $H$ with $H(0)=0$.
 Pakes and Navarro (2007) show under some mild conditions that the distribution function $H_{\alpha,\beta}$
 of $X$ determines $H$. Based on that result we derive in this paper
 a recursive formula for calculation of $H$, if $H_{\alpha,\beta}$ is known. Furthermore,
  we investigate the relation between the tail asymptotic behaviour of $X$ and $Y$.
We present three applications of our asymptotic results concerning the extremes of
two random samples with underlying distribution functions $H$ and $\HAB$, respectively,
and the conditional limiting distribution of  bivariate elliptical distributions.

{\it Key words and phrases:} Beta random scaling; fractional
integral; elliptical distribution; max-domain of attraction;
asymptotics of sample maxima; conditional limiting results,
estimation of conditional distribution; Weibull-tail distribution;
Gardes-Girard estimator.

\section{Introduction}
Let $X,Y,B$ be three independent random variables such that
\BQN\label{1} X &\equaldis & BY , \EQN where $\equaldis$ stands for
equality of the distribution functions. In our context the random
variable $B$ plays the role of a random scaling or multiplier.  Clearly, if the distribution functions of $Y$ and $B$ are
known, then the distribution function of $X$ can be easily
determined. In various theoretical and practical situations the
question of interest is whether the distribution function of $Y$ can
be determined provided that those of $X$ and $B$ are known. Indeed,
random scaling of $Y$ by $B$ is treated in several papers and
different contexts, see for instance the recent contributions
Tang and Tsitsiashvili (2003,2004), Jessen  and  Mikosch(2006), Tang (2006,2008), Pakes (2007), Pakes and Navarro (2007),    Beutner and Kamps (2008a,b).

Unless otherwise stated, in this article we fix $B$ to be a beta random variable with
 positive parameters $\alpha,\beta$. If $H$ denotes the distribution
function of $Y$, then the distribution function of $X$ (denoted by $H_{\alpha,\beta})$ is defined in
terms of $H$ and both parameters $\alpha, \beta$. If $Y$ is another beta random
variable, then $X$ is the product of two such beta random variables, which have  been studied extensively in
the literature, see Galambos and Simonelli (2004), Nadarajah (2005),
Nadarajah and  Kotz (2005b, 2006), Dufresne (2007), Beutner and Kamps (2008a) and the references therein.

Our main impetus for dealing with the beta random scaling comes from
 Pakes and Navarro (2007) which paves the
way for the distributional and asymptotic considerations in this
paper. Theorem 2.2 therein gives an explicit formula for the
calculation of the distribution function $H$, provided that $\HAB$
satisfies some weak growth restrictions on its derivatives. Utilising the
aforementioned theorem, we show in this paper that the distribution
function $H$ can be calculated iteratively without imposing any
additional assumption on $\HAB$. This iterative inversion may lack the elegance of the explicit formula
in Pakes and Navarro (2007),
but it turns out to be quite useful in asymptotic contexts
where we can define the tail behaviour of the survivor
function of $Y$ when that of the
survivor function of $X$ is known, and vice-versa.

We present three applications of our asymptotic results:\\
a) Determining which maximal domain of attraction contains $H_{\alpha,\beta}$ when the membership of $H$ is known;\\
b) The derivation of conditional limiting results for bivariate elliptical random vectors; and\\
c) New estimators for the conditional
distribution function and the conditional quantile function
of bivariate elliptical random vectors allowing one component of the random vector to grow to infinity.

The paper is organized as follows.  In the next section we give some
preliminary results. The main result of Section 3 is the iterative
inversion for $\HAB$ --  \netheo{eq:mainth} below. In Section 4 we
investigate the asymptotic relation of the survivor function of $X$ and $Y$ under conditions arising in
extreme value theory, showing in particular that $H$ is attracted to an extreme value distribution if and only if $\HAB$ is attracted to the same distribution. The direct implications are formulated (in Section 7) in a generality which subsumes the particular case of beta scaling. Conditional limiting results and estimation of
conditional distribution function for bivariate elliptical random vectors is discussed in Sections 5 and 6. All proofs  and some related results are relegated to Section 7.

\section{Preliminaries}
We introduce notation and then discuss some properties of
the Weyl fractional-order integral operator. A
key result of Pakes and Navarro (2007) is recalled because it is  crucial for our
considerations.

We use notation such as $X \sim F$ to mean that $X$ is a random variable with distribution function $F$, and $\overline F:=1- F$ denotes  the corresponding  survivor function. The upper endpoint of the distribution
function $F$ is denoted by $r_F$ and its lower endpoint by $l_F$.
If  $\alpha,\beta >0$  then
$\mathrm{beta}(\alpha,\beta)$ and  gamma$(\alpha,\beta)$ denote respectively the beta and the gamma
distributions  with corresponding density
functions
$$(B(\alpha,\beta))^{-1}x^{\alpha-1}(1- x)^{\beta-1}, \quad x\in (0,1),\quad
\text{ and   } \frac{\beta^{\alpha}}{\Gamma(\alpha)}x^{\alpha-1}
\exp(- \beta x) , \quad x\in (0,\IF),$$ where $B(\alpha,\beta)$ is the beta function and 
$\Gamma(\alpha)$ is the gamma function. Since beta distributed random variables appear below
in several instances, we use exclusively the notation $B_{\alpha,\beta}$ for
a beta random variable with parameters $\alpha,\beta$. On occasion it is convenient to extend the definition to understand $\pk{B_{0,\beta}=0}=1$ if $\beta>0$ and $\pk{B_{\alpha,1}=1}=1$ if $\alpha>0$. Unless otherwise stated, factors in products of random variables are assumed to be independent.\\

Next, define the Weyl fractional-order integral operator
$I_\beta, \beta
>0 $ by
\BQN (I_\beta h) (x)&:=&\frac{1}{\Gamma(\beta)}\int_x^\infty
(y-x)^{\beta-1} h(y) \, dy, \quad x>0,
\EQN
with ${h}:[0,\IF) \to \R$ a measurable function. The function $I_\beta h$  is well
defined  if (see Pakes and Navarro (2007))
$$ \int_{\ve}^\IF x^{ \beta- 1} \abs{h(x)}\, dx < \IF$$
is satisfied for all $\ve >0$, in which case we write $h\in
\kal{I}_\beta$ with the understanding that $\beta$ may assume
negative values. Define further (consistently) $I_0 h :=h$. If $h$
is a density function of a positive random variable $Y \sim H$,
then $I_\beta h$ is well-defined for every $\beta>0$. Suppose  $g$ is a measurable function such that if $Y\sim H$, then
$\E{Y^{\beta-1}|g(Y)|}<\infty$. Then we define
\BQN   (\OPJ_{\beta, g} H)(x) &=&\frac{1}{\Gamma(\beta)}\int_x^{r_H}(y-x)^{\beta-1}g(y) \, d
H(y) , \quad \forall x\in (l_H,r_H), \EQN
i.e.,  $\OPJ_{\beta,g} $ denotes the Weyl-Stieltjes
fractional-order integral operator acting on the class of
distribution functions on $\R$ with weight function $g$.

The Weyl fractional-order integral operator is closely
related to beta random scaling. To see this, let $ \alpha,\beta>0$ and $Y>0$ and $\bab$ be independent random
variables such that
\BQN \label{eq:1}
X & :\equaldis & Y B_{\alpha,\beta}, \quad \text{ where } X \sim \HAB, Y \sim H,
\EQN
and $l_H\ge 0$. In the light of equation (14) in Pakes and Navarro
(2007), for any $x\in (l_H,r_H)$ we have
\BQN \label{PKS:14}
\HAB(x)&=&\frac{\Gamma(\alpha+ \beta)}{\Gamma(\alpha)} x^{\alpha}
(I_\beta p_{- \alpha -\beta} H)(x),
\EQN
with $p_{s}$ the power function defined by
$$ p_{s}(x):= x^{s},\quad s\inr, \ x>0.$$
We mention in passing two important topics in probability theory and
statistical applications where the Weyl fractional-order
integral operator is encountered: a) the sized- or length biased law
(see e.g., Pakes (2007), Pakes and Navarro (2007)); and  b) the Wicksell
problem (see e.g., Reiss and Thomas (2007)).
For the essentials of fractional integrals and derivatives see Miller and Ross (1993).

Now we state three properties of{$I_\beta h$.}
\BL \label{lem:a0}
Let $\beta,c$ be positive constants, and let $h$ be a real measurable function. \\
i) If $h\in \kal{I}_{\beta+ c}$, then
\BQN
I_\beta  I_c h =   I_c  I_\beta h = I_{\beta +c} h.
\EQN
ii) Let  $D^n$ denote the $n$-fold derivative operator ($n\inn$). If the $n$-fold derivative $h^{(n)}:= D^n h$ exists
almost everywhere and $h^{(n)} \in \kal{I}_\beta$, then
\BQN\label{DN}
   D^n I_\beta h &=& I_\beta h^{(n)}
\EQN and
\BQN D^k I_{n} &=& (-1)^k I_{n-k}, \quad k=1 \ldot n. \EQN
iii) If  $\lambda \in (0,\beta)$ and $H$  is a distribution function on $\R$ with $H(0)=0$, then
\BQN \label{eq:IIa}
( I_{\beta- \lambda}   p_{- \beta}   (I_\lambda p_{-\alpha- \lambda} H) )(x) &=&
 x^{- \lambda} (I_\beta p_{-\alpha- \beta} H)(x), \quad \forall x \in (0, \infty).
\EQN  \EL
The next theorem, which is an insignificant
variation of Theorem 2.2 in Pakes and Navarro (2007) shows that the
survivor function $\overline H$ can be retrieved by applying the
differential and the Weyl fractional-order integral
operator to $\BARHAB$.

\BT \label{theo:0} Let $H$,$\HAB, \alpha,\beta \in (0,\IF)$ be as above, with  $\HAB(0)=0$.
If  $\HAB^{(n-1)}$ is absolutely continuous and $\HAB^{(n-i)} \in \kal{I}_{\delta- \alpha- i}, i=0 \ldot n$ with
$\delta$ and $n$ such that
\BQN \beta + \delta =:n\inn, \quad \delta \in [0,1),
\EQN
 then
\BQN\label{pakes}
\overline H(x)&=& (-1)^n \frac{\Gamma(\alpha)}{\Gamma(\alpha+\beta)}
x^{\alpha+ \beta} (I_\delta D^n p_{-\alpha} \BARHAB)(x)
\EQN
holds for any $x\in (0, r_H).$
\ET
\section{Iterative Calculation of $H$}
Let $X,Y, B_{\alpha, \beta}$, related by  \eqref{eq:1},  be as above. In this section our main interest is the determination of
 $H$ from the known form of  $\HAB$. As already mentioned, an explicit formula is presented as
Theorem 2.2 in Pakes and Navarro (2007) (see \eqref{pakes} above).
If $\beta \in (0,1]$, then the only requirement for the validity of their
theorem is that $\HAB(0)=0$, which obviously is fulfilled whenever $H(0)=0$.
The following well-known multiplicative property of beta random variables is the key   to our iterative version of  Theorem 2.2 above. Specifically, if $\lambda \in (0,\beta)$, then
$$ B_{\alpha,\beta} \equaldis B_{\alpha,\lambda} B_{\alpha+ \lambda, \beta- \lambda}.$$
 Consequently, \eqref{eq:1} implies that
 \BQN
\label{eq:2} X \equaldis YB_{\alpha, \beta} & \equaldis & Y
B_{\alpha,\lambda} B_{\alpha+ \lambda, \beta - \lambda}. \EQN
Theorem 2.2 of Pakes and Navarro (2007) and \eqref{eq:2}implies
the following result:

\BT \label{lem:a1} Let $\alpha,\beta$ be two positive constants, and let $X,Y, B_{\alpha,\beta}$ be
independent random variables satisfying \eqref{eq:1} with $X\sim \HAB, Y\sim H$, and $H(0)=0$. \\
i) If $\lambda \in (0,\beta),$ then
\BQN \label{eq:thm:1:A}
\BARHAB(x)
 &=&  \frac{\Gamma(\alpha+ \beta)}{\Gamma(\alpha)} x^{\alpha+\lambda}
 ( I_{\beta- \lambda}   p_{- \beta}
 (I_\lambda p_{-\alpha- \lambda} \overline H  ))(x),
 \quad \forall x\in (0,r_H).
\EQN \\
ii) If $\beta- \lambda \in [0,1),$ and $\delta\in [0,1)$ is such that $\beta- \lambda+ \delta= 1$, then
\BQN
\label{eq:thm:1:B}
\BARHAL (x) &=& \frac{\Gamma(\alpha+\lambda)}{\Gamma(\alpha+ \beta)} x^{\alpha+ \beta}
 \Bigl[ (\alpha+ \lambda) (I_\delta p_{-\alpha- \lambda-1} \BARHAB)(x)+
 (\OPJ_{\delta,  p_{-\alpha- \lambda} } \HAB)(x) \Bigr],
 \quad \forall x\in (0, r_H).
\EQN \ET We state next a simple corollary which is of some interest
in the context of the Weyl fractional-order integral
operator.
\BK \label{kor:A} Let $H$ be a distribution function on
$\R$ such that $H(0)=0$. Then for any $x\in (0,r_H)$ we have
\BQN
x^{\alpha-1} (\OPJ_{\beta, p_{-\alpha -\beta+1}}H )(x) &=& (D
p_{\alpha} I_\beta p_{- \alpha -\beta} H) (x) = - (D
p_{\alpha} I_\beta p_{- \alpha -\beta} \overline H) (x).
\EQN Moreover, if $H$ possesses the density function $h$, then
\BQN
(I_\beta p_{- \alpha -\beta+1}h )(x)
&=&  \alpha  (I_\beta
p_{- \alpha -\beta} H) (x)+ x (I_\beta D( p_{- \alpha -\beta}
H)) (x), \quad x\in (0,r_H).
\EQN
\EK
\COM{ \BQN
   x^{\alpha-1} (I_\beta p_{- \alpha -\beta-}h )(x) &=&
    D(x^{\alpha} (I_\beta p_{- \alpha -\beta} H)) (x) =
    \alpha x^{\alpha -1}
    (I_\beta p_{- \alpha -\beta} H)) (x)+ x^{\alpha}     (I_\beta D( p_{- \alpha -\beta} H)) (x).
\EQN }
The main result of this section is the following iterative
formula for computing  $H$ when  $\HAB$  is known.

\BT \label{eq:mainth}
Let $X \sim \HAB, Y \sim H$ and $B_{\alpha,\beta}, \alpha,\beta>0$ be three independent random
variables satisfying \eqref{eq:1} such that $\HAB(0)=0$. If $\beta_0:=\beta> \beta_1 >\cdots >
\beta_k> \beta_{k+1}:=0,$ with $k\in \{0,\N\}$ and $ \delta_i,i\le k+1$ are constants such that
\BQN
\lambda_i:= \beta_{i-1}- \beta_{i} \in (0,1],  \quad \delta_i:=1- \lambda_i, \quad i=1 \ldot k+1,
\EQN
then we can construct distribution functions $H_0:=H, H_1 \ldot H_{k+1}=H_{\alpha,\beta}$ such that
\BQN
 \overline H_{i-1}(x) &= & \frac{\Gamma( \alpha + \beta_{i} )}{\Gamma(\alpha + \beta_{i-1})}
 x^{\alpha + \beta_{i-1}} \Bigl[ (\alpha+ \beta_{i}) (I_{\delta_{i}} p_{- \alpha - \beta_{i}-1} \overline H_{i})(x)
      +(\OPJ_{\delta_{i},  p_{- \alpha - \beta_{i}}} H_{i} )(x)\Bigr], \quad \forall x\in (0,r_H).
 \EQN
\ET

\begin{remark} \label{rem:2}
(a) Let $ B_i {\sim B_{\alpha_i,\beta_i}, i\ge 1}$ be independent beta random variables
and  independent of $Y \sim H$. If the random variable $X$ with distribution function $H_n, n\ge 2$
has the stochastic representation
\BQN
 X \equaldis Y \prod_{i=1}^n B_i^{c_i}, \quad c_i\in (0,\IF), \quad  i=1 \ldot n,
 \EQN
then  \netheo{eq:mainth}  implies that  $H$ can be retrieved recursively
from  $H_n$, provided that $H_n(0)=0$.

(b) An interesting (open) question arises in connection with random products. Specifically,
if $\kal{N}$ is a counting random variable taking positive integer values independent of $Y, B_i,i\ge 1,$
such that
\BQN\label{eq:rN}
 X& \equaldis &Y \prod_{i=1}^{\kal{N}} B_i^{c_i} \quad \text{ where } X \sim H_{\kal{N}},
 \EQN
then under what conditions on $\kal{N}$ can we (recursively) compute the distribution function $H$ if $H_{\kal{N}}$ is known?\\
{Also arises a similar  question if $X,Y$ are related by
\BQN\label{eq:rNa}X &\equaldis &Y [ B_3 B_1+ (1- B_3)B_2].
\EQN
}
%
\end{remark}

\section{Tail Asymptotics}
The tail asymptotics of products have been studied in papers such as
Berman (1983, 1992), Cline and Samorodnitsky (1994), Tang and Tsitsiashvili (2003, 2004), Jessen  and  Mikosch (2006), Tang (2006, 2008), and the references therein. Our asymptotic considerations below can be motivated
by considering sample maxima.

Specifically, let $X_i, Y_i, i=1 \ldot n,$ be independent copies of $X=Y B_{\alpha,\beta}$ and $Y$, respectively, and
$$M_{X,k}:= \max_{1 \le j \le k} X_j,  \quad M_{Y,k}:= \max_{1 \le j \le k}Y_j, \quad k\ge 1 $$
be the corresponding sample maxima. From extreme value theory (see e.g., de Haan and Ferreira (2006), Falk et al.\
(2004, p. 23), Resnick (1987, p. 38)) if there are constants $a_n>0, b_n$ such that
\BQN\label{eq:maxH}
 \limit{n} \sup_{t\inr} \Abs{H^n(a_nt+ b_n)- Q(t)} &=&0,
\EQN
then we have the convergence in distribution
\BQN (M_{Y,n}-
b_n)/a_n &\todis& \kal{M}_{Y} \sim Q, \quad n\to \IF,
\EQN
where $Q$ is a univariate extreme value distribution (Gumbel, \FRE or
Weibull). If $\eqref{eq:maxH}$ holds (write $H \in MDA(Q)$) it is of
some interest to investigate the asymptotic behaviour of
$M_{X,k},k\ge 1$, where $X_i,i \le n$ are the results of a beta
random scaling i.e.,
\BQN\label{eq:i} X_i &\equaldis &Y_i
B_i, \quad B_i \equaldis B_{\alpha,\beta}, \quad i=0 \ldot n, \quad
n \ge 1,
\EQN
with  $Y_i\sim H$ and the $B_i$'a and $Y_i$'s mutually independent. Thus $X_i\sim \HAB$.  A key question is whether $\HAB$ is in a maximal domain of attraction if $H$ is, and conversely? We answer this below, as well as  exposing the explicit tail asymptotic relations underlying \eqref{eq:i}.
\subsection{Gumbel Max-domain of Attraction}
If \eqref{eq:maxH} holds with $Q= \Lambda$
the unit Gumbel distribution ($\Lambda(x):=\exp(-\exp(-x)),x\inr$),
then there exists a positive measurable {\em scaling} function $w$ (see e.g., de Haan and Ferreira (2006), Resnick (1987, p. 46))
 such that
 \BQN \label{eq:gumbel}
 \lim_{x \uparrow r_H} \frac{\overline H(x+t/w(x))} {\overline H(x)} &=& \exp(-t),\quad \forall t\inr
\EQN
is valid.  We write $H\in MDA(\Lambda, w) $ if \eqref{eq:gumbel} holds. The scaling function $w$ satisfies
\BQN\label{eq:uv:B}
 \lim_{x\uparrow r_H} x w(x)=\infty, \quad \text{and  }
\lim_{x\uparrow r_H} w(x)(r_H - x) =\infty, \quad \text{if  } r_H < \infty,
\EQN
and also the self-neglecting property
\BQN\label{eq:self}
\lim_{x\uparrow r_H}  \frac{w(x+ t/w(x))}{w(x)} = 1,
\EQN
which holds locally uniformly for $t\in \R$; see e.g., Resnick (1987, p.\ 41). Note that most authors work with the so-called auxiliary function $1/w(x)$, but our convention follows Berman (1992) because results we prove are closely linked to some in  his Chapter 12.

Canonical examples of distribution functions in the Gumbel max-domain of attraction are the univariate Gaussian and the gamma distributions, which are special cases of distribution functions whose scaling functions have the form (for $x$ large)
\BQN\label{eq:Abd} w(x)& =&\frac{r \theta
x^{\theta-1}}{1+ L_1(x)},
\EQN
where  $L_1(x)$ is  regularly varying at infinity with index $\theta \mu, \mu \in (-\IF, 0)$
and $r,\theta$ are positive constants. Note that $\theta=2$ for the Gaussian case, and we have for the $\mathrm{gamma}(\alpha,\beta)$ case that $\theta=1$,  $w(x)=\beta $ and
\BQN \label{eq:gumbel:2}
 \limit{x} \frac{\overline H(x+t)} {\overline H(x)} &=& \exp(- \beta t),\quad \forall t\inr.
\EQN
Distribution functions $H$  that satisfy \eqref{eq:gumbel:2} comprise  what in other contexts is called the exponential tail class $\kal{L}(\beta)$.  See Pakes (2004) for references, and Pakes and Steutel (1997) where they are called medium-tailed.

We state now the first result of this section, a close relative of Theorem 12.3.1 in Berman (1992); see Example 1 below for the latter.  In \S7 we will state and prove the general proposition \netheo{theo:new:G} which subsumes both direct assertions.
\BT \label{eq:theo:BM1:0}
Let $H, \HAB$ be as in \netheo{eq:mainth}. Then $H\in MDA(\Lambda,w) $  iff (if and only if) $\HAB\in MDA(\Lambda,w) $. If one of these holds, then
\BQN
\label{res:BERM1:a:01}
\BARHAB(x) &=& (1+o(1))K (xw(x))^{-\overlineAL} \overline H(x) , \quad x \uparrow r_H,
\EQN
where $K:=\Gamma(\overlineALI+\overlineAL)/\Gamma(\overlineALI) $, and the density function $\HABD$  of $\HAB$ satisfies
\BQN \label{eq:Misses:1}
\lim_{x\uparrow r_H} \frac{\HABD(x)}{w(x)\BARHAB(x)} &=&1.
\EQN
\ET

The asymptotic equivalence \eqref{res:BERM1:a:01} is the principal assertion here, as can be seen by noting that if one of the distribution functions $F$ and $H$ is in $MDA(\Lambda,w)$ and they are related by
\BQN \label{eq:gumrel}
\o F(x)=(1+o(1))x^c(w(x))^\mu \o H(x),\qquad(x\uparrow r_H),
\EQN
where $c,\mu$ are real, then it follows from \eqref{eq:gumbel} and \eqref{eq:self} that the other distribution function is in $MDA(\Lambda,w)$.

{It is well-known that if $H$ is a univariate distribution function with upper endpoint $r_H=\IF$ and
$H\in MDA(\Lambda, w)$, then $\overline{H}$ is rapidly varying (see Resnick (1987)) i.e.,
\BQN\label{eq:rapid}
\limit{x} \frac{\overline{H}(cx)}{\overline{H}(x)}&=&0,  \quad \forall c> 1.
\EQN
A necessary ingredient in the proof of  \netheo{eq:theo:BM1:0} is the following rate of convergence refinement to \eqref{eq:rapid}; recall the first member of \eqref{eq:uv:B}.
\BL\label{rapid} Let $H$ be a univariate distribution function with $r_H=\IF$. If $H\in MDA(\Lambda, w)$,
then we have for any constant $\mu\ge 0$
\BQN\label{eq:rapid:mu}
\limit{x} (x w(x))^\mu  {\overline{H}(cx) \over  \overline{H}(x)} &=&0,  \quad \forall c> 1.
\EQN
\EL
}

\begin{remark}

(a) The self-neglecting property \eqref{eq:self} implies that the density function $\HABD$ of $\HAB$ satisfies
\BQNY
 \frac{\HABD(x+ t/w(x))}{\HABD(x)}& \to & \exp(-t), \quad x\uparrow r_H
\EQNY
locally uniformly for $t\in \R$, provided that either $H\in MDA(\Lambda, w)$, or $\HAB \in MDA(\Lambda, w)$.

(b) By \netheo{eq:theo:BM1:0},  if $\HAB\in MDA(\Lambda, w),$ then we can reverse \eqref{res:BERM1:a:01} obtaining
\BQN
\label{res:BERM1:a:02}
\overline H(x) &=& (1+o(1))\frac{\Gamma(\overlineALI)}{\Gamma(\overlineALI+\overlineAL)}
(xw(x)) ^{\overlineAL}\BARHAB  (x)  ,
\quad x \uparrow r_H.
\EQN
{See Berman (1992)  and Hashorva (2007d) for similar results.
Further note that \eqref{eq:rapid:mu} and \eqref{res:BERM1:a:02} imply for any $c\in (1,\IF)$ that
$$\overline H(x)= o(\BARHAB  (cx)), \  \ \BARHAB  (x)=o(\overline H (cx)),  \text{  and} \ \ (xw(x)) ^{\overlineAL}\BARHAB  (x)=o(1), \quad x\uparrow r_H.
$$}
\end{remark}

We give next two illustrations  of \netheo{theo:new:G}.\\
{\bf Example 1.} (a) Theorem 12.3.1 in Berman (1992) follows from \netheo{theo:new:G}(a) by taking (see  \eqref{eq:betaTail})
$$\phi(u)=\pk{\sqrt{1-B_{\alpha,\beta}}>u}$$
 and checking that, since $1-\bab\eqdis B_{\beta,\alpha}$,  \eqref{eq:betaTail} holds with $C=2^\alpha/\alpha B(\alpha,\beta)$ {\it and} the exponent $\beta$ replaced with $\alpha$.

(b) Let $H, F$ be two distribution functions as in \netheo{eq:theo:BM1:0} and suppose that
 $l_H=0$ and $r_H=\IF$. We assume that the random multiplier $B$ has the
stochastic representation
\BQNY
B  \equaldis \lambda U_1+ (1- \lambda)U_2, \quad \lambda \in (0,1),
\EQNY
where $U_1,U_2$ are two independent positive random variables such that for $i=1,2$
$$ \pk{U_i> 1- s}= (1+o(1))c_i s^{d_i}, \quad  c_i,d_i \in (0,\IF), \quad s\downarrow 0.$$
It follows that as $s\downarrow 0$
\BQNY
 \pk{B  > 1- s} &=& (1+o(1))\frac{c_1 c_2}{\lambda^{d_1} (1- \lambda)^{d_2}} \frac{\Gamma(d_1+1)\Gamma(d_2+1)}{ \Gamma(d_1+d_2+1)}s^{d_1+d_2} .
\EQNY
Further, assume for all large $x$ that
\BQN\label{eq:kotz:Fk}
\overline{H}(x)&=&  (1+o(1))M  x^{N}\exp(-r
x^\theta ), \quad M >0,r>0,\theta >0, N\inr.
 \EQN
Since, for any $t\inr$,
\BQNY
\limit{x}\frac{ \overline{H}( x+ t x ^{1-\theta }/(r \theta )   )}{\overline{H}(x)}&=& \exp(-t)
\EQNY
we have $H\in MDA(\Lambda, w)$ with
\BQN\label{eq:sc:III}
w(x)&=& r \theta  x ^{\theta -1},\quad x>0.
\EQN
In view of \netheo{theo:new:G} the distribution function $F$ of $BY$ satisfies $F \in MDA(\Lambda,w)$ and, as $x\to \IF$,
\BQNY
\o F  (x)&=&    (1+o(1))   C^*x^{N-\theta(d_1+ d_2)} \exp(-r x^\theta ),
 \EQNY
with
\BQNY
 C^*&=&  M  (r \theta)^{-d_1- d_2} \frac{c_1 c_2}{\lambda^{d_1} (1- \lambda)^{d_2}} \Gamma(d_1+1)\Gamma(d_2+1).
 \EQNY

\COM{
 is the density function $\HABD$ of $\HAB$ has tail asymptotics
\BQNY\HABD(x)&=&  (1+o(1))Kr \theta    x^{N+\theta -1} \exp(-rx^\theta ), \quad x\to \IF.
 \EQNY
Furthermore, also $H$ is in the Gumbel max-domain of attraction with the same scaling function $w$, and
\eqref{res:BERM1:a:02} implies 
\BQNY
\overline H(x) 
&=&  (1+ o(1)) K (r \theta) ^{\overlineAL} \frac{\gamma(\overlineALI)}{\gamma(\overlineALI+\overlineAL)}
x ^{N+ \overlineAL \theta}     \exp(-r x^\theta ),  \quad x \to \IF .
\EQNY
}

\subsection{Regularly Varying Tails}
We deal next with distribution functions $\HAB$ in either the \FRE or the Weibull max-{domains}
of attraction.  As we will discuss below, the asymptotics of $\HAB$ when $H$ is attracted to the \FRE distribution is quite well known, and results for the  Weibull max-domain of attraction  are less complete. In Section 7 we offer  simpler proofs of these results, and their converses, i.e., when $\HAB$ belongs to one of these max-domains of attractions, then so does $H$.

The unit \FRE distribution function with positive index $\gamma$ is $\Phi_\gamma(x):=\exp(-x^{-\gamma}),x>0 $. It is well-known that a distribution function $H$ with infinite upper endpoint $r_H=\IF$
is in the \FRE max-domain of attraction (see e.g., Falk et al.\ (2004), Resnick (1987)) iff $\overline H$
is regularly varying {at infinity} with index $-\gamma<0$, i.e.,
\BQN\label{eq:F:gamma}
\limit{x} \frac{\overline H(xt)}{\overline H(x)}= t^{-\gamma}, \quad \forall t\in (0,\IF).
\EQN
If $l_H=0$ and $0<\gamma<1$, then this condition is the criterion that $H$ is attracted to a positive stable law with index $\gamma$. Breiman (1965, Proposition 3) shows that if this holds, then the distribution function $F$ of $X=BY$, where the random multiplier $B$ is independent of $Y$, is also attracted to the same positive stable law provided that  $\E{|B|}<\infty$. (Thus $B$ is not restricted in sign or magnitude.) Specifically,  $H$ and $F$ are tail equivalent, i.e.,
\BQN\label{eq:FreTE}
\o F(x)= (1+o(1)) \E{B^\gamma} \o H(x) \quad(x\to\infty).
\EQN
Jessen and Mikosch (2006, p.\ 184) observe that Breiman's proof is valid for any positive $\gamma$ if $B\ge0$ and  $\E{B^{\gamma+\epsilon}}<\infty$ for some $\epsilon>0$. Berman (1992, Theorem 12.3.2) proves this tail equivalence for the case $B=\sqrt{1-\bab}$.

So in particular, we conclude that if  $ \alpha, \beta>0$ and  $\HAB,$ is  defined via \eqref{eq:1} with $\HAB(0)=0$, then
\BQN \label{res:BERM1:b2}
\BARHAB(x)&=& (1+o(1))\E{B_{\alpha,\beta}^\gamma}\overline H(x), \quad x\to \infty,
 \EQN
and
$$\E{B_{\alpha,\beta}^\gamma}= \frac{\Gamma(\barALI+\barAL)\Gamma(\barALI+\gamma)}
{\Gamma(\barALI)\Gamma(\barALI+\barAL+\gamma)}.$$

The next theorem asserts that this tail equivalence holds also if $\gamma=0$, and conversely,  if  $\gamma>0$ and $\HAB \in MDA(\Phi_\gamma)$, then so is  $H$. Breiman's methodology is completely analytical, and in Section 7 we shall give  a  much simpler proof for the case of a general bounded multiplier  $0\le B\le 1$. We indicate too how this can be extended to the general result.
\BT \label{eq:theo:BM1:2} Let $H,\HAB, \alpha, \beta>0$ be two distribution functions defined via
\eqref{eq:1} with $H(0)=0$. Then $H$ satisfies \eqref{eq:F:gamma} with some $\gamma\ge 0,$ iff
$\HAB$  satisfies \eqref{eq:F:gamma} with the same index $\gamma$.
Furthermore, for any $\gamma >0$ we have
\BQN\label{eq:Misses:2}
\limit{x}  \frac{x h_{ \alpha,\beta} (x)}{\BARHAB(x)} &=& \gamma.
\EQN
\ET

%
%

{\bf Example 2.} \netheo{eq:theo:BM1:2} shows in particular that Pareto tails are preserved under independent beta random scaling.

The unit Weibull distribution function with index $\gamma>0$ is  $\Psi_\gamma(x):= \exp(-\abs{x}^\gamma), x< 0$.  It is well known that if  $H$ has a finite upper endpoint (say $r_H=1$), then  $H\in MDA(\Psi_\gamma)$ iff
\BQN\label{eq:weib:gamma}
\lim_{x \downarrow 0} \frac{\overline H(1-tx)}{\overline H(1- x)}
&=& t^{\gamma} , \quad \forall t>0.
 \EQN

Theorem 12.3.3 in Berman (1992) is closely related to the following result, and in Section 7 we prove a general theorem which subsumes both direct assertions.
\BT \label{eq:theo:BM1:3} Let $H, \HAB, \alpha,\beta$ be as in \netheo{eq:theo:BM1:2}. If $r_H=1, H(0)=0$
and \eqref{eq:weib:gamma} holds for some $\gamma \ge 0$, then $\HAB{ \in MDA( \Psi_{\overlineAL+\gamma})}$ and
\BQN
\label{res:BERM3} \BARHAB(1- x)&=&
(1+o(1))K
x^{\overlineAL} \overline H(1-x), \quad x \downarrow 0,
\EQN
with $K:= \Gamma(\overlineALI+\overlineAL) \Gamma(\gamma+1)/(\Gamma(\overlineALI)\Gamma(\gamma+\overlineAL+1)$.\\
Furthermore we have
\BQN \label{res:BERM3:d}
\lim_{x \downarrow 0}\frac{x h_{\alpha,\beta}(1-x)}{\BARHAB(1- x) }&=& \overlineAL+\gamma>0.
\EQN
Conversely, if {$\HAB \in MDA(\Psi_{\overlineAL+\gamma}),\gamma\ge 0$}, then \eqref{eq:weib:gamma} is satisfied.
\ET

\bigskip

\begin{remark}
(a) If \eqref{eq:1} holds with $B_{\alpha,\beta}\sim \mathrm{gamma}(\alpha,\beta)$,
then in Lemma 17 of  Hashorva et al.\ (2007) it is shown that $\overline H$ satisfies \eqref{eq:F:gamma} with some  $\gamma\ge 0,$  iff  $\BARHAB$ satisfies \eqref{eq:F:gamma} with the same index $\gamma$
 {(see also Jessen  and  Mikosch (2006)).}

(b) Under the Gumbel or the Weibull max-domain of attraction assumption on  $H$ or $\HAB$ by \eqref{eq:uv:B} we have
$$ \lim_{x \uparrow r_H}\frac{\BARHAB  (x)}{\overline H(x)} =0,$$
whereas when $H$ or $\HAB$ are in the \FRE max-domain of attraction the above limit is a positive constant.
\end{remark}

\COM {\bibitem{fang}
 (1990) ...  In {\it Statistical inference in
elliptically contoured and related distributions},  K.T. Fang and T.W. Anderson, eds,
Allerton Press, New York, pp.\ 127--136.
}

\COM{
{\bf Example 3.} Let $\HAB$ be a beta distribution function with positive parameters $\lambda, \delta$, which depend
on $\alpha,\beta$, and denote by $h_{\alpha,\beta}$ its density function.

Assume that $H$ relates to $\HAB$ via \eqref{eq:2}  satisfies \eqref{eq:weib:gamma} with $\gamma\ge 0$. By \eqref{res:BERM3:d}
we have
\BQNY
\lim_{ x\downarrow 0}\frac{x h_{\alpha,\beta}(1-x)}{\BARHAB (1- x)}&=& \delta= \beta+ \gamma,
\EQNY
hence $\gamma= \beta- \delta$. Further we have
\BQNY
\label{res:BERM3}\overline H(1-x)
\COM{
 &=&(1+o(1))\frac{\gamma(\barALI)\gamma(\gamma+\barAL+1)}{\gamma(\barALI+\barAL) \gamma(\gamma+1)}
 x^{-\barAL} \BARHAB(1- x)\\
&=&(1+o(1))\frac{\gamma(\barALI)\gamma(\gamma+\barAL+1)}{\gamma(\barALI+\barAL) \gamma(\gamma+1)}
\frac{\gamma(\lambda+\delta)}{\gamma(\lambda)\gamma(\delta+1)}   x^{\beta+ \gamma -\barAL} \\
&=&(1+o(1))\frac{\gamma(\barALI)\gamma(\gamma+\barAL+1)}{\gamma(\barALI+\barAL) \gamma(\gamma+1)}
\frac{\gamma(\lambda+\beta + \gamma)}{\gamma(\lambda)\gamma(\beta + \gamma+1)}   x^{\gamma } \\
&=&(1+o(1))\frac{\gamma(\barALI)}{\gamma(\barALI+\barAL) \gamma(\gamma+1)}
\frac{\gamma(\lambda+\beta + \gamma)}{\gamma(\lambda)}   x^{\gamma } \\
}
&=&(1+o(1))\frac{\gamma(\barALI)\gamma(\lambda+\beta + \gamma)}{\gamma(\lambda)\gamma(\barALI+\barAL) \gamma(\gamma+1)}
   x^{\gamma }, \quad x \downarrow 0,
\EQNY
hence $\gamma$ is necessarily strictly positive. 
In the special case $H=beta(\alpha_1,\beta_1 )$, we have further $\beta_1=\gamma= \delta- \beta$ and
\BQNY
\label{res:BERM3}\overline H(1-x)
&=&(1+o(1))\frac{\gamma(\barALI)\gamma(\lambda+\delta)}{\gamma(\lambda)\gamma(\barALI+\barAL) \gamma(\delta- \beta+1)}
   x^{\gamma }, \quad x \downarrow 0
\EQNY
implying
$$ \frac{\gamma(\alpha_1+ \delta- \beta)}{\gamma(\alpha_1)} =
\frac{\gamma(\barALI)\gamma(\lambda+\delta)}{\gamma(\lambda)\gamma(\barALI+\barAL) } . $$
Consequently, $\alpha_1:=\alpha+\beta $ is an admissible parameters, provided that $\alpha=\lambda$.
}

\section{Conditional Limiting Results}
Let the bivariate random vector $(O_1,O_2)$ be uniformly distributed on the unit circle, $R\sim H$ be independent of $(O_1,O_2)$, and let $(S_1,S_2) \equaldis R (O_1,O_2)$ be the corresponding bivariate (planar) spherical random vector. Finally, define the bivariate elliptical random vector
 \BQN \label{e}
 (U,V) \equaldis  (S_1, \rho S_1+ \sqrt{1- \rho^2} S_2), \quad \rho \in (-1,1).
 \EQN

Distributional properties of spherical and elliptical random vectors are studied by many authors, e.g., Cambanis et al.\ (1981), Fang et al.\ (1990), Kotz et al.\ (2000) and their references.  Referring to Cambanis et al.\ (1981) we have
 \BQN\label{eee}
 O_1^2 &\equaldis& O_2^2 \sim \mathrm{beta}(1/2,1/2).
 \EQN
Basic asymptotic properties of spherical and elliptical random vectors can be derived utilising \eqref{e} and \eqref{eee}.
One line of enquiry  is to determine the  asymptotic behaviour of the conditional distribution of  $V-\rho U$ given an event constraining the values of $U$. For example,  in several statistical applications (see Abdous et al.\ (2005)) the approximation of the conditional random variable
$$Z_x^*\equaldis  (V - \rho x) \lvert U > x,\quad x\inr$$
 is of some interest.  Since $V-\rho U= \sqrt{1- \rho^2} S_2$, the outcome follows directly from Theorem 12.3.3 in Berman (1992), i.e., if $H\in MDA(\Lambda,w)$, then
   \BQN\label{ajm:2}
c(x) Z_x^* & \todis & \sqrt{1- \rho^2} Z \quad x\uparrow  r_H,
\EQN
where $c(x):= \sqrt{w(x)/x},x>0,$ and $Z$ is a standard Gaussian random variable.   Abdous et al. (2005) is an independent account.
 \netheo{eq:ellip:A} below embellishes this outcome.

The point-wise conditioned random variable
$$Z_x\equaldis (V - \rho x) \lvert U = x,\quad x\inr$$
is a particular case of the conditional  multivariate models introduced by Heffernan and Tawn  (2004) for treating certain inference problems.  They raise the general problem of conditional limit laws when one component of a random vector tends to infinity, and they give results for some particular parametric families. It is known that (Hashorva (2006), Corollary 3.1) that $Z_x$ has the same Gaussian limit law as $Z_x^*$, i.e.,
\BQN\label{ajm:1}
 c(x) Z_x \todis \sqrt{1- \rho^2} Z, \quad x\uparrow  r_H.
\EQN
We will prove that if $H$ is absolutely continuous then \eqref{ajm:1} holds in the stronger sense that the density functions converge. We prove in addition that both limit assertions hold assuming that the (marginal) distribution of $|U|$ is attracted to the Gumbel distribution. Finally, Hashorva and Kotz (2009) gives an account of these results based on the strong Kotz approximation.

\BT \label{eq:ellip:A} Let $ H,(U,V), \rho\in
(-1,1),c(x),Z_x,Z_x^*, x>0$ be as above with $|U| \sim G$ and $G(0)=0$. If
$G\in MDA(\Lambda, w)$ or $H\in MDA(\Lambda, w)$, then (a),  \eqref{ajm:2} is satisfied; and (b), \eqref{ajm:1} is satisfied if, in addition, $H$ is absolutely continuous.
\ET

The proof of this theorem rests on a closure lemma for distributions attracted to the Gumbel distribution.

\BL\label{power} Let $0\le X\sim F$, $p>0$ be a constant, and denote the distribution function of $X^p$ by $F_p$. Then $F\in MDA(\Lambda,w)$ iff $F_p\in MDA(\Lambda,w_p)$ where
$$w_p(x)=p^{-1} x^{(1/p)-1}w\left(x^{1/p}\right). $$
\EL


\section{Estimation of Conditional Survivor and Quantile Function}
\def\cre#1{#1}
For $i=1,2,\dots$, let $(U_i,V_i)$ be independent copies of $(U,V)$ as defined in the previous section, and suppose too that $R\sim H\in MDA(\Lambda,w)$ with $r_H=\infty$.   We are interested in the conditional survivor function
\BQNY
{\Psi_x}(y)&:= &\pk{V> y \lvert U> x}, \quad x,y\inr.
\EQNY
 Estimation of the distribution function $1- \Psi_x(y)$
when $x$ is large is discussed in detail by Abdous et al.\
(2007). As noted there, if $x$ is large  there may be insufficient data
available for the effective  estimation of $\Psi_x(y)$. Similar difficulties apply
for estimation of the inverse function (or conditional quantile function), $\Theta(x,\cdot),s\in (0,1),x\inr$ of $1- \Psi_x(\cdot)$. The Gaussian approximation implied by \netheo{eq:ellip:A} entails
\BQN\label{gaussA} \sup_{y\inr} \Abs{\Psi_x( y \crA{\sqrt{x/w(x)}}+
\rho x) - \Phi( \crA{y} /\sqrt{1- \rho^2})} \to  0, \quad x\to
\IF, \EQN
where $\Phi$ is the standard Gaussian distribution function.

On this basis, Abdous et al. (2007) propose two estimators of $\Psi_x$.  \netheo{eq:ellip:A} implies that the Gaussian approximation in
\eqref{gaussA} is valid if we assume instead that $U\sim G\in MDA(\Lambda,w)$. For estimation purposes this fact is
crucial because we can estimate $w$ based only on the random sample $U_1 \ldot U_n,$ or $V_1 \ldot V_n$.

A non-parametric estimator  $\hat \rho_n$ of $\rho$ is given by (see e.g., Li and Peng (2009)))
\BQN\label{eq:Bl}
 \hat \rho_n&:= &\sin(\pi \hat \tau_n /2), \quad n>1,
 \EQN
where $\hat \tau_n$ is the empirical estimator of Kendall's tau.  Now, if $\hat w_n(x)$ is an
estimator of the scaling function $w(x)$ (for all large $x$), then by the above
approximation we can estimate $\Psi_x(y)$ by
\BQN\label{eq:estim}
\hat \Psi_{n,x}(y)&:= &\overline \Phi \Bigl(\hat h_n (y- \hat \rho_n x )/(1- \hat \rho_n^2)^{1/2}
\Bigr), \quad n> 1,
\EQN
where $\hat h_n(x):=(\hat w_n(x)/x)^{1/2}, x>0.$
 An estimator for the quantile function $\Theta$ is then given by
\BQN
 \hat \Theta_n(x,s) &=&  \rho_n x+  \sqrt{1- \hat \rho_n^2}\Phi^{-1}(s)/\hat h_n(x), \quad x>0, s\in (0,1),
 \EQN
with $\Phi^{-1}$ the inverse of $\Phi$. Both of these estimators are consequences of the Gaussian approximation. However, our concern here is with  estimation of $w$. Specifically, we assume that the scaling function $w$ satisfies \eqref{eq:Abd} with positive constants $r,\theta$ and $L_1$  regularly varying  with
index $\theta \mu, \mu \in (-\IF, 0)$. It follows that  (see Abdous et al.\ (2007))
\BQN
\label{eq:KotzM} \o G(x) &=& \exp(- r x^\theta(1+ L_2(x))
\EQN
holds for all large $x$, where $L_2$ is another regularly varying function  with
index $-\theta \mu.$ This places $G$ in the class of Weibull-tail
distributions, and $\theta^{-1}$ is the so-called Weibull tail-coefficient (see Gardes and Girard (2006), or
Diebolt et al.\ (2007)). Canonical examples of  Weibull-tail distributions are the
Gaussian, gamma,  and extended Weibull {distributions.} Next, define for $i=1 \ldot n$,
$$R_i^{(1)}:= U_{i},\quad  R_{in}^{(2)}:=
\sqrt{ U_{i}^2+ (V_{i} - \hat \rho_n U_{i})^2/(1- \hat \rho_n^2)}$$
and write $R_{1:n}^{(k)}\le \cdots \le R_{n:n}^{(k)}, k=1,2$ for the associated order statistics. Based on $R_i^{(1)}, i\le n$ or $R_{in}^{(2)}, i\le n$ we may construct
the Gardes-Girard (2006)  estimator of $\theta$,
$$ \hat \theta_n^{(j)}:=
\frac{1}{T_n}\frac{1}{k_n}\sum_{i=1}^n \Bigl( \log R_{n-i+1:n}^{(j)}- \log R_{n- k_n+1:n}^{(j)}\Bigr), \quad j=1,2,$$
where $1 \le k_n \le n, T_n>0, n\ge 1$ are constants satisfying
$$ \limit{n} k_n=\IF, \quad \limit{n} \frac{k_n}{n}=0, \quad
\limit{n} \log(T_n/k_n)= 1, \quad \limit{n} \sqrt{k_n} b(\log(n/k_n))\to \lambda\inr,$$
and the function $b$ (related to $L_1$) is  regularly varying with index $\eta$.  The scaling coefficient $r$ can be estimated by (see Abdous et al.\ (2007))
\BQN
\hat r_{n}^{(j)}&=& \frac{1}{k_n} \sum_{i=1}^{k_n}
\frac{ \log (n/i) }{(R_{n-i+1:n}^{(j)})^{\hat \theta_n^{(j)}}}, \quad j=1,2, n>1,
\EQN
leading to the following estimators of  $w$,
\BQN
\hat w_{n}^{(j)}(x)&=& \hat r_n ^{(j)} \hat \theta_n ^{(j)} x^{\hat \theta_n ^{(j)}- 1}, \quad x>0, j=1,2, n>1.
\EQN
Our suggestion is to estimate $w$ by $\hat w_n^{(1)},$ because it is based on independent and identically distributed $R_{i},i\le n$. This differs from the estimator $\hat w_n^{(2)}$ recommended by Abdous et al.\ (2007) which is based on the  dependent random variables  $R_{1n}^{(2)}\ldot  R_{nn}^{(2)}$ (recall $\hat \rho_n$
is estimated from $(U_i,V_i),i\ge 1$).

A third estimator of  $w$ can be easily constructed by considering the sample $V_1 \ldot V_n$ since
by the assumption $U\equaldis V$.

Note in passing that if $\theta=1$, then we have the estimator of $r$ (of interest for $G$ in $\kal{L}(r),r>0$)
\BQN
\hat r_{n}^{(1)}&=& \frac{1}{k_n} \sum_{i=1}^{k_n}
\frac{ \log (n/i) }{R_{n-i+1:n}^{(1)}}, \quad  n>1.
\EQN

\section{Further Results and Proofs}
We present first some asymptotic results for the Weyl
fractional-order integral  operator, followed by the proofs of all the results in the previous sections.

\BT \label{theo:OPJ:G} Let $H$ be a univariate distribution function with $H(0)=0,$ $r_H\in (0,\IF]$, and
$H\in MDA(\Lambda,w)$. If $\alpha$ is real and $\beta>0$,  then
\BQN\label{eq:opj}
(\OPJ_{\beta, p_{-\alpha}} H)(x)&=& (1+o(1)) (w(x))^{-(\beta -1)}x^{- \alpha} \overline H(x) , \quad
x\uparrow r_H,
\EQN and
\BQN\label{eq:opj:2} (I_\beta p_{- \alpha} \overline H)(x)
&=& \frac{(1+o(1))}{ w(x)} (\OPJ_{\beta,p_{-\alpha}} H)(x) , \quad x\uparrow r_H.
\EQN
\ET
\prooftheo{theo:OPJ:G} Let $W_x$ be a random variable whose survivor function is
$$\pk{W_x>z}={\o H(x+z/w(x)) \over  \o H(x)}, \qquad (0\le z<r(x)),$$
 where
$$
r(x)=(r_H-x)w(x) \ \mathrm{if} \ r_H<\infty, \quad\&\quad  =\infty \  \mathrm{if} \ r_H=\infty.
$$
Then \eqref{eq:gumbel} is equivalent to the convergence assertion $W_x\todis W$  which has the standard exponential distribution. Observe now that  for $x\in (0,r_H)$ we may write
\BQNY
(\OPJ_{\beta, p_{-\alpha}}H) (x)&=& \frac{1}{\Gamma(\beta)}\int_x^{r_H} (y-x)^{\beta -1} x^{- \alpha}\, d H(y)\\
&=& {x^{-\alpha}\o H(x) \over \Gamma(\beta)}\int_0^{r(x)} \left({z\over w(x)}\right)^{\beta-1}
(1+z/v(x))^{-\alpha} d_z\pk{W_x\le z},
\EQNY
where $v(x)=xw(x)$ and  we have used the substitution $y=x+z/w(x)$ for the second equality. Hence
$$ {(w(x))^{\beta-1} x^\alpha \over \o H(x)} (\OPJ_{\beta, p_{-\alpha}}H) =
{1\over \Gamma(\beta)}\vek{E}\left\{ W_x^{\beta-1}\left(1+W_x/v(x)\right)^{-\alpha}\right\}.$$
It follows from \eqref{eq:uv:B} and the moment convergence theorem (Feller (1971, p. 252)) that the expectation converges to $\E{W^{\beta-1}}=\Gamma(\beta)$. This proves \eqref{eq:opj}.

The same manoeuvres yield
$${x^\alpha(w(x))^{\beta} \over \o H(x)} (I_\beta p_{- \alpha} \overline H)(x)=
{1\over \Gamma(\beta)}\int_0^{r(x)} z^{\beta-1} (1+z/v(x)^{-\alpha} \pk{W_x>z} dz\to 1,$$
using the dominated convergence theorem, and \eqref{eq:opj:2}  follows.  \QED

\netheo{theo:OPJ:G} subsumes and generalizes results in Berman (1992, \S12.2) applying to the case $r_H=\infty$. To align  with Berman's notation, we use $\beta-1$ to denote his parameter $p$, and in what follows we assume that
$\E{Y^\beta}<\infty$.\\
(i) Propositions 12.2.3 and 4 in Berman (1992) concern distribution functions $F$ having the form
$$\o F(x)=(1+o(1))c\int_x^\infty (y-x)^{\beta-1} \o H(y)dy.$$
It is easily seen that
$$\o F(x)=(1+o(1))c\Gamma(\beta)(\OPJ_{\beta+1,p_0}H)(x)=(1+o(1))c\Gamma(\beta)(w(x))^{-\beta} \o H(x),$$
and this is valid if $\beta>0$, which extends the range of parameter in Berman's Proposition 12.2.4. \\
(ii) Proposition 12.2.5 in Berman (1992) concerns survivor functions proportional to the order-$q$ stationary excess distribution generated by $H$, i.e.,
$$\o F(x)=(1+o(1))c\int_x^\infty y^{q-1}\o H(y)dy,$$
where $q$ is real. The integral can be recast as
$$q^{-1}\int_x^\infty \left(y^q-x^q\right) dH(y)=q^{-1} x^q\o H(x)\E{(1+W_x/v(x))^q-1},$$
from which it follows that, as $x\to\infty$,
$$\o F(x)=(1+o(1))cx^{q-1}{\o H(x) \over w(x)}.$$
(iii) The order-$q$ size-biased distribution generated by $H$ induces survivor functions of the form
$$\o F(x)=(1+o(1)) c\int_x^\infty y^q dH(y)=(1+o(1))c(\OPJ_{1,p_q}H)(x)
=(1+o(1))c x^q\o H(x).$$
It follows from \eqref{eq:gumrel} that each  above $F\in MDA(\Lambda,w)$.

Note that if $q>0$, then the results under (ii) and (iii) are related {\em via} \netheo{eq:theo:BM1:0} because if $\widehat Y_q$ and $\widetilde Y_q$ denote the order-$q$ size-biased and stationary excess versions of $Y$, then $\widetilde Y_{q}\eqdis \widehat Y_q B_{q,1}$. See Pakes (1996, \S4) for this connection and further generalization involving beta scaling.\\

%
%
\BT \label{lem:fre:1} Let $H$ be a univariate distribution function
with  $r_H=\IF$. Assume that $H(0)=0$ and \eqref{eq:F:gamma} holds with $\gamma\ge 0$.
If $\beta>0$ and $c$ are two constants such that $\beta+ c< \gamma+1$, then
\BQN\label{eq:fre:lem:1}
(\OPJ_{\beta, p_{c}} H)( x ) &=& (1+o(1)) \frac{\gamma \Gamma( \gamma+1-
\beta - c)}{\Gamma( \gamma+1-c)}\overline H(x)x^{\beta+c-1}, \quad x \to\IF.
 \EQN
Furthermore if $\gamma \ge 0$, then
\BQN\label{eq:fre:lem:2}
(I_{\beta}  p_{c} \overline H)( x)
&=& \frac{\Gamma(\gamma-
\beta-c)}{\Gamma(\gamma-c)} \overline H(x)x^{\beta+c}  , \quad x \to \IF.
\EQN
\ET
\prooftheo{lem:fre:1} Let $W_x$ have the distribution function $\max(1-\o H(xt)/\o H(x),0)$. Then \eqref{eq:F:gamma} is equivalent to: If $\gamma>0$, then $W_x\todis W$ which has the Pareto survivor function $t^{-\gamma}$ for $t\ge 1$; and if $\gamma=0$, then $W_x\toprob \infty$.

Substituting $y=tx$ into the integral defining $\OPJ_{\beta,p_c}H$ gives the representation
$$(\OPJ_{\beta,p_c}H)(x)={x^{\beta+c-1} \o H(x) \over \Gamma(\beta)}\vek{E}\left[(W_x-1)^{\beta-1}W_x^c\right].$$
If $\gamma>0$ and $\epsilon>0$ is chosen so $\beta+c+\epsilon<\gamma+1$, then $\vek{E}(W^{\beta+c+\epsilon-1})<\infty$, and hence the above expectation converges to
$$\vek{E}\left[(W-1)^{\beta-1}W^c\right]=\gamma B(\gamma+1-\beta-c,\beta),$$
and \eqref{eq:fre:lem:1} follows. This assertion follows too if $\gamma=0$ because $(W_x-1)^{\beta-1}W_x^c<W_x^{\beta+c-1}$, and the exponent is negative.

The same substitution yields
$$(I_{\beta,p_c}\o H)(x)={\o H(x) x^{\beta+c} \over \Gamma(\beta)}\int_1^\infty (t-1)^{\beta-1} t^c \pk{W_x>t}dt,$$
and it is clear that the integral converges to
$$\int_1^\infty (t-1)^{\beta-1}t^{c-\gamma}dt=B(\gamma-\beta-c,\beta),$$
whence \eqref{eq:fre:lem:2}.\QED\\

%
\BT \label{lem:weib:1} Let $H$ be a univariate distribution function
with upper endpoint $r_H=1$.  {Assume that $H(0)=0$,  and that \eqref{eq:weib:gamma} holds with
 $\gamma \ge 0$}.  If  $ \beta>0$ and $ c\inr$ are constants and $\gamma >0$, then
\BQN\label{eq:weib:lem:1}
(\OPJ_{\beta, p_{c}} H)( 1- x) &=&(1+o(1))\frac{\Gamma(\gamma+1)}{\Gamma(\beta+ \gamma)} \overline H(1- x)
x^{\beta-1} , \quad x \downarrow 0
\EQN
and  if $\gamma\ge 0$, then
\BQN\label{eq:weib:lem:2}
(I_{\beta}  p_{c} \overline H)( 1- x) &=&(1+o(1))\frac{1}{\gamma+ \beta} x  (\OPJ_{\beta, p_{c}} H)(1- x) ,
\quad x \downarrow 0.
\EQN
\ET
\prooftheo{lem:weib:1} Let $W_x\le1$ be a random variable having the distribution function $\o H(1-tx)/\o H(1-x)$. If $\gamma>0$, then \eqref{eq:weib:gamma} is equivalent to $W_x\todis W:=U^{1/\gamma}$, where $U$ has the standard uniform distribution (i.e. beta$(1,1)$), and if $\gamma=0$, then $W_x\todis 1$.  The substitution $y=1-xt$ yields
$$(\OPJ_{\beta, p_{c}} H)( 1- x)  = {\o H(1-x) x^{\beta-1} \over \Gamma(\beta)}\E{(1-W_x)^{\beta-1}(1-xW_x)^{c}}.$$
If $\gamma>0$, then the expectation converges as $x\downarrow$ to
$$\vk{E}\left\{(1-W)^{\beta-1}\right\}= \gamma B(\gamma,\beta),$$
and if $\gamma=0$ then it converges to unity.  So \eqref{eq:weib:lem:1} follows in both cases.

The same substitution yields
$$ (I_{\beta}  p_{c} \overline H)( 1- x) ={\o H(1-x) x^{\beta} \over \Gamma(\beta)}\int_0^1 (1-t)^{\beta-1}(1-xt)^{c}
\pk{W_x\le t}dt,$$
and the integral converges to $B(\gamma+1,\beta)$.  Thus
$$(I_{\beta}  p_{c} \overline H)( 1- x)= (1+o(1))\overline H(1- x) x^{\beta} \frac{\Gamma(\gamma+1)}{\Gamma(\beta+ \gamma+1)},$$
and \eqref{eq:weib:lem:2} follows. \QED
%

\bigskip
\prooflem{lem:a0}
Since the first two statements are borrowed from Lemma 2.2 in Pakes and Navarro (2007)
we show next only statement $iii)$. Let $Y\sim H, B_{\alpha, \beta} , B_{\alpha,\lambda}$ and $B_{\alpha+ \lambda, \beta - \lambda}$ be independent random variables. For any $\lambda \in (0,\beta)$
we have the stochastic representation (see \eqref{eq:2})
\BQNY
Y B_{\alpha, \beta} & \equaldis & Y^* B_{\alpha+ \lambda, \beta -
\lambda}, \quad Y^* \equaldis Y B_{\alpha,\lambda},
\EQNY
with $Y^*\sim H_{\alpha, \lambda}$ another random variable independent of
$B_{\alpha+ \lambda, \beta - \lambda}$. Applying \eqref{PKS:14} we obtain for any $x\in (0, r_H)$
\BQNY \HAB(x)&=& \frac{\Gamma(\alpha+ \beta)}{\Gamma(\alpha)} x^{\alpha} (I_\beta p_{- \alpha -\beta} H)(x)\\
&=& \frac{\Gamma(\alpha+ \beta)}{\Gamma(\alpha+ \lambda)} x^{\alpha+
\lambda} (I_{\beta- \lambda} p_{- \alpha -\beta} \HAL)(x)
\EQNY
and
\BQNY
\HAL (x)&=& \frac{\Gamma(\alpha+ \lambda)}{\Gamma(\alpha)}
x^{\alpha} (I_\lambda p_{- \alpha -\lambda} H)(x).
\EQNY
Consequently
\BQNY
 (I_\beta p_{- \alpha -\beta} H)(x) &=&
x^{\lambda} (I_{\beta- \lambda} p_{ -\beta}  (I_\lambda p_{-
\alpha -\lambda} H))(x),
 \EQNY
and the result follows. \QED
\COM{ \BQNY F(x)&=& \frac{\gamma(\alpha+ \beta)}{\gamma(\alpha+
\lambda)} x^{\alpha+ \lambda}
(I_{\beta - \lambda} \IDo_{- \alpha -\beta} H^*)(x)\\
&=& \frac{\gamma(\alpha+ \beta)}{\gamma(\alpha+ \lambda)} x^{\alpha+
\lambda} \frac{\gamma(\alpha+ \lambda)}{\gamma(\alpha)}
(I_{\beta - \lambda} \IDo_{- \alpha -\beta} \IDo_{\alpha} (I_\lambda\IDo_{- \alpha -\lambda} H) )(x)\\
&=& \frac{\gamma(\alpha+ \beta)}{\gamma(\alpha)} x^{\alpha+ \lambda}
(I_{\beta - \lambda}\IDo_{-\beta} (I_\lambda \IDo_{- \alpha -\lambda} H) )(x)\\
&=& \frac{\gamma(\alpha+ \beta)}{\gamma(\alpha)} x^{\alpha} (I_\beta
\IDo_{- \alpha -\beta} H)(x). \EQNY \BQNY x^{\lambda} (I_{\beta -
\lambda} \IDo_{-\beta} (I_\lambda\IDo_{- \alpha -\lambda} H) )(x)&=&
(I_\beta \IDo_{- \alpha -\beta} H)(x). \EQNY }

\prooftheo{theo:0} The proof follows immediately from Theorem 2.2 in
Pakes and Navarro (2007) and the identity
\BQN \label{eq:IIDO}
(I_{\beta,p_{-c}})(x)={\Gamma(c-\beta) \over \Gamma(c)} x^{\beta-c}
\EQN
  \QED


\prooftheo{lem:a1} The identity \eqref{eq:IIDO} implies that
\BQNY
 1 &= &\frac{\Gamma(\alpha+ \beta)}{\Gamma(\alpha)} x^{\alpha} (I_\beta p_{- \alpha -\beta})(x),
 \quad \forall x\in (0, r_H),
 \EQNY
and hence \eqref{PKS:14}
\BQNY
\BARHAB(x)&=& 1- H_{\alpha,\beta}(x)\\
&=& \frac{\Gamma(\alpha+ \beta)}{\Gamma(\alpha)} x^{\alpha} (I_\beta
p_{- \alpha -\beta})(x)-
\frac{\Gamma(\alpha+ \beta)}{\Gamma(\alpha)} x^{\alpha} (I_\beta p_{- \alpha -\beta} H)(x)\\
&=& \frac{\Gamma(\alpha+ \beta)}{\Gamma(\alpha)} x^{\alpha} \Bigl[
(I_\beta p_{- \alpha -\beta})(x)- (I_\beta p_{- \alpha -\beta} H)(x)\Bigr]\\
&=& \frac{\Gamma(\alpha+ \beta)}{\Gamma(\alpha)} x^{\alpha} (I_\beta
p_{- \alpha -\beta} \overline H)(x), \qquad (x\in (0, r_H)) \EQNY
thus the  first result follows utilising further \eqref{eq:IIa} which holds if $\overline H$ replaces $H$.

We show next the second claim. Since $H(0)=0$,  Lemma 2.1 in Pakes and Navarro (2007) shows that
$$\HAL(0)= \HAB(0)=0.$$
Furthermore, both $\HAL$ and $\HAB$ are absolutely continuous and
$$ X \equaldis Y^* B_{\alpha+ \lambda, \beta- \lambda}, \quad \text{ with  } Y^*\sim H_{\alpha,\lambda}, \quad X \sim \HAB.$$
 Therefore, in order to show the proof we need to check the assumptions of \netheo{theo:0}. In our case $n=1$, hence the
condition $\HAB^{(n-1)}=\HAB^{(0)}=\HAB$ is absolutely continuous is satisfied.
Since $\HAB^{(1)}$ is a density function and $\delta\in [0,1),$ then clearly $ \HAB^{(1)}\in \kal{I}_{\delta- \alpha- \lambda }$. Further we have $ \HAB^{(0)}=  \HAB \in \kal{I}_{\delta- \alpha- \lambda-1 }$ since $\HAB$ is bounded by 1. Applying \netheo{theo:0} for any $x\in (0, r_H)$ we may write
\BQNY \BARHAL(x)&=& - \frac{\Gamma(\alpha+ \lambda)}{\Gamma(\alpha+ \beta)}
 x^{\alpha+ \beta} (I_\delta D (p_{-\alpha- \lambda } \BARHAB))(x)\\
 &=& \frac{\Gamma(\alpha+ \lambda)}{\Gamma(\alpha+ \beta)} x^{\alpha+ \beta}
 \Bigl[ (\alpha+ \lambda)(I_\delta p_{-\alpha- \lambda-1} \BARHAB)(x)+
(\OPJ_{\delta,  p_{-\alpha- \lambda} } \HAB)(x) \Bigr],
\EQNY and the result follows.
\QED

\proofkorr{kor:A} Letting $\lambda \to 0$ in \eqref{eq:thm:1:A} we obtain (recall  {$I_0 h:= h$})
\BQN\label{eq:91} \BARHAB(x) &=& \frac{\Gamma(\alpha+ \beta)}{\Gamma(\alpha)} x^{\alpha} (I_\beta p_{- \alpha -\beta} \overline H)(x), \quad
 \forall x\in (0,r_H).
\EQN
Consequently, we have
\BQNY
-\HABD(x)&=& \frac{\Gamma(\alpha+ \beta)}{\Gamma(\alpha)}
(D (p_{\alpha} I_\beta p_{- \alpha -\beta} \overline H))(x), \quad
 \forall x\in (0,r_H)
 \EQNY
and in view of \eqref{PKS:14},
\BQNY\HABD(x)&=& \frac{\Gamma(\alpha+ \beta)}{\Gamma(\alpha)} (D(
p_{\alpha} I_\beta p_{- \alpha -\beta} H))(x) , \quad  \forall x\in (0,r_H).
\EQNY
Since $\HABD$ is given by (see (22) in Hashorva et al.\ (2007))
\BQN \label{PKS:14:2}
\HABD(x)
&=& \frac{\Gamma(\alpha+ \beta)}{\Gamma(\alpha)} x^{\alpha-1}  (\OPJ_{\beta ,p_{- \alpha -\beta+1}} H)(x), \quad  \forall x\in (0,r_H),
\EQN
the result follows. \QED

\prooftheo{eq:mainth}
\COM{$$ Y_k B(\alpha, \beta_{k}) = Y$$
$$ Y_k = \prod_{i=1}^k  B(\alpha+  \beta_i, \lambda_i) Y
\BQN H(x)&=& \frac{\gamma(\alpha+ \beta_1)}{\gamma(\alpha+ \beta)}
x^{\alpha+ \beta}
 \Bigl[ (\alpha+ \lambda)I_\delta (\IDo_{-\alpha- \lambda-1} \BARHAB)(x)+
(\OPJ_{\delta,  \IDo_{-\alpha- \lambda} } \HAB)(x) \Bigr], \quad \forall x\in (0, r_H)\\
\EQN } Let $B_{\alpha+\beta_i, \lambda_i}\sim \mathrm{beta}(\alpha+\beta_i,
\lambda_i), i=0 \ldot k$
be independent beta random variables  independent of $X$ and $Y$. 
By the assumptions we may write \BQNY
  X &\equaldis &Y B_{\alpha, \beta_0}\\
  &\equaldis&
  B_{\alpha, \beta_1}  Y B_{\alpha+  \beta_1, \beta_0- \beta_1}\\
  &\equaldis &    Y_1 B_{\alpha, \beta_1}, \quad \text{ with }
  Y_1 \equaldis  Y_0 B_{\alpha+  \beta_1, \beta_0- \beta_1} \equaldis  Y_0 B_{\alpha+  \beta_1, \lambda_1}, \quad Y_0:=Y.
  \EQNY
Similarly \BQNY
  X &\equaldis& Y_2    B_{\alpha, \beta_2}, \quad  \text{ with }Y_2 \equaldis  Y_1 B_{\alpha+  \beta_2, \lambda_2}
  \EQNY
and repeating we arrive at \BQNY
  X &\equaldis& Y_k    B_{\alpha, \beta_k}, \quad  \text{ with } Y_k \equaldis  Y_{k-1} B_{\alpha+  \beta_{k}, \lambda_{k}}.
  \EQNY
Setting $Y_{k+1}:=X$ we may write the above stochastic
representation as \BQNY
 Y_{k+1} &\equaldis& Y_k    B_{\alpha+ \beta_{k+1}, \lambda_{k+1}}.
 \EQNY
Let $H_0:=H$ and $ H_{k+1}:= H_{\alpha,\beta}$. Applying \eqref{eq:thm:1:B} we obtain for any $i=1 \ldot k+1$,
\BQN
 \overline H_{i-1}(x) &= \frac{\gamma( \alpha + \beta_{i} )}{\gamma(\alpha + \beta_{i}  +\lambda_{i})}
 x^{\alpha + \beta_{i-1}} \Bigl[ (\alpha+ \beta_{i}) (I_{\delta_{i}} p_{- \alpha - \beta_{i}-1} \overline H_{i})(x)
 +(\OPJ_{\delta_{i},  p_{- \alpha - \beta_{i}}} H_{i} )(x)\Bigr],
 \EQN
and the assertion follows. \QED

We precede our account of scaling relations for the Gumbel distribution with the following proof.\\
\prooflem{rapid} If $\beta>0$ it follows from \netheo{theo:OPJ:G} that
$$\o H_{1,\beta}(x)=(\OPJ_{\beta+1,p_\beta})(x)= (1+o(1)) \Gamma(1+\beta) {\o H(x) \over (xw(x))^\beta},\qquad(x\uparrow r_H).$$
On the other hand, if $B:=B_{1,\beta}$ and $c>1$, then
$$\o H_{1,\beta}(x)>\int_{cx}^\infty  \pk{B>x/y} dH(y)> \pk{B>c^{-1}} \o H(cx).$$
Combining these estimates yields
$$\limsup_{x\to\infty}(xw(x))^\beta {\o H(cx) \over \o H(x)}<\infty.$$
The assertion \eqref{eq:rapid:mu} follows by choosing $\beta>\mu$ and appealing to \eqref{eq:uv:B} in the case $r_H=\infty$.\QED

The next result is the foreshadowed  generalization of the direct assertion of \netheo{eq:theo:BM1:0}. It comprises two parts which respectively yields a tail estimate of the distribution function of a random scaling, and its density function.
\BT \label{theo:new:G} Suppose $H\in MDA(\Lambda,w)$. (a) If $\phi(u)\ge 0$ is defined and bounded on $[0,1]$ and it satisfies
\BQN\label{eq:betaTail}
\phi(u)= (1+o(1)) C(1-u)^\beta, \quad u\uparrow 1,
\EQN
where $\beta,C\ge 0$ are constants, then
$$I(x):=\int_x^\infty \phi(x/y)dH(y)= (1+o(1))C\Gamma(1+\beta){\o H(x) \over (xw(x))^\beta},\quad x\uparrow r_H.$$
(b) If $g(u)\ge 0$ is defined in $[0,1]$ such that $ug(u)$ is defined and bounded on $[0,u']$ for any $u'<1$, and
\BQN\label{eq:betaTail:b}
g(u)=(1+o(1))c(1-u)^{\beta-1}\quad u\uparrow 1,
\EQN
where $c\ge0$ and $\beta>0$ are constants, then
$$J(x):= \int_x^\infty y^{-1}g(x/y)dH(y)=(1+o(1))c\Gamma(\beta){\o H(x) \over x^\beta (w(x))^{\beta-1}},\quad  x\uparrow r_H.$$
\ET
\prooftheo{theo:new:G} We prove only (b) since the details for (a) are similar and simpler.  If $u'\in(0,1)$, then $x/y\le u'$ if $y\ge x/u'$ and
$$J_1(x):=\int_{x/u'}^\infty y^{-1} g(x/y)dH(y)=O[x^{-1} \o H(x/u')].$$
If $r_H$ is finite, then $J_1(x)=0$ if $x>u'r_H$. If $r_H=\infty$, then, recalling that $v(x)=xw(x)$,  \nelem{rapid} ensures that $J_1(x)=o\left(x^{-1} \o H(x)(v(x))^{-\mu}\right)$ ($x\to\infty$) for all positive $\mu$.

If $c$ is positive and $0<\epsilon\ll c$, then it follows from \eqref{eq:betaTail:b} that $g(u)/(1-u)^{\beta-1}\in (c-\epsilon,c+\epsilon)$ if $u'<u<1$ and $u'$ is sufficiently close to unity. Hence  $J(x)-J_1(x)$ is asymptotically equal to
$$J_2(x):=c\int_x^{x/u'} y^{-1}(1-x/y)^{\beta-1} dH(y).$$
Proceeding as in the proof of \netheo{theo:OPJ:G} we obtain the representation
$$J_2(x)= {c \o H(x) \over x(v(x))^{\beta -1}}  \vk{E}\left\{{W_x^{\beta-1} \over (1+W_x/v(x))^\beta}; W_x\le v(x)(1-u')/u'\right\}.$$
The expectation converges to $\E{W^{\beta-1}}=\Gamma(\beta)$. Taking $\mu>\beta$ above, we see that $J_1(x)=o(J_2(x))$, and the assertion follows.\QED

\prooftheo{eq:theo:BM1:0} Assume that $H\in MDA(\Lambda,w)$.  The direct assertion \eqref{res:BERM1:a:01} follow from \netheo{theo:new:G}(a) by setting $\phi(u):=\pk{\bab>u}$ and checking that \eqref{eq:betaTail} holds with $C=[\beta B(\alpha,\beta)]^{-1}$. Next, taking $g(u)$ as the density function of $\bab$ it is obvious that the conditions of  \netheo{theo:new:G}(b) are satisfied with $c=1/B(\alpha,\beta)$. Thus \eqref{eq:Misses:1} follows from \eqref{res:BERM1:a:01} and \eqref{eq:betaTail:b}.

To prove the converse, assume that $\HAB\in MDA(\Lambda,w)$  for some positive scaling function $w$. With
the notation of \netheo{eq:mainth} we may write for $i=1 \ldot k+1$
\BQN\label{eq:hik}
 \overline H_{i-1}(x) &= & \frac{\Gamma( \alpha + \beta_{i} )}{\Gamma(\alpha + \beta_{i-1})}
 x^{\alpha + \beta_{i-1}} \Bigl[ (\alpha+ \beta_{i}) (I_{\delta_{i}} p_{- \alpha - \beta_{i}-1} \overline H_{i})(x)
      +(\OPJ_{\delta_{i},  p_{- \alpha - \beta_{i}}} H_{i} )(x)\Bigr], \quad \forall x\in (0,r_H),
 \EQN
where $\overline H_0:=\overline H, \overline H_{k+1}:= \BARHAB $. In view of
\netheo{theo:OPJ:G} and \eqref{eq:uv:B}, we obtain for
$i=k+1$ that
 \BQNY
 \overline H_{i-1}(x) &= & (1+o(1))\frac{\Gamma( \alpha + \beta_{i} )}{\Gamma(\alpha + \beta_{i-1})}
 x^{\alpha + \beta_{i-1}} (\OPJ_{\delta_{i},  p_{- \alpha - \beta_{i}}} H_{i} )(x)\\
&= & (1+o(1))\frac{\Gamma( \alpha + \beta_{i} )}{\Gamma(\alpha +
\beta_{i-1})}
 x^{\beta_{i-1}- \beta_{i}} (w(x))^{-(\delta_i-1)}\overline H_{i} (x) , \quad \forall x \uparrow r_H.
 \EQNY
By \eqref{eq:uv:B} and \eqref{eq:self} it follows that $H_{k}\in MDA(\Lambda,w)$.
Since the above holds for all $i=1 \ldot k$, it follows that
$H_0=H \in MDA(\Lambda,w)$ too. Next,  \eqref{eq:91} and \eqref{PKS:14:2} imply for any  $x>0$ that
 \BQN \label{ratio:h:H}
 \frac{  h_{\alpha,\beta}(x)}{ \BARHAB(x)}&=&
 \frac{x(\OPJ_{\beta, p_{- \alpha- \beta}} H)(x)}{(I_\beta p_{- \alpha- \beta -1} H)(x) },
 \EQN
so applying \netheo{theo:OPJ:G} establishes \eqref{eq:Misses:1}, and the result follows.
\QED


As foreshadowed above, the following argument  includes a simple proof of \eqref{eq:FreTE}  for an arbitrary bounded random scaling. We then show how this proof can be extended to remove the boundedness restriction.\\
\prooftheo{eq:theo:BM1:2} With $W_x$ as in the proof of \netheo{lem:fre:1}, clearly
$$\pk{XB>x}=\o H(x) \E{\pk{B>W_x^{-1}}}.$$
But
$$\pk{B>W_x^{-1}}\to \pk{B>W^{-1}}=\int_0^1 \pk{W>u^{-1}}d\pk{B\le u}=\E{B^\gamma},$$
and hence $\o F(x)=(1+o(1))\o H(x)\E{B^\gamma}$.

Similarly, if $B$ has the density function $g(u)$, then the density function of $F$ is
$$f(x)=x^{-1} \o H(x) \E{W_x^{-1}g(W_x^{-1})}.$$
If $g$ satisfies appropriate boundedness conditions, which certainly are satisfied by beta density functions, then the expectation converges to
$$\vek{E}\left\{W^{-1}g(W^{-1})\right\}=\gamma \int_1^\infty w^{-\gamma-2}g(w^{-1})dw=\gamma \E{B^\gamma}.$$
It follows that $(h(x)/\o H(x))=(1+o(1))(\gamma/x)$. The direct assertions of \netheo{eq:theo:BM1:2} follow.

The converse asserts that if $\BARHAB$ is regularly varying with index $-\gamma\le0$, then $\overline H$ is also regularly varying with index $-\gamma$. If $\gamma >0,$ then the proof follows from \netheo{eq:mainth} and \netheo{lem:fre:1}.

Alternatively, write  can write $\bab=Z_1/(Z_1+Z_2)$, where $ Y, Z_1,Z_2, $ are independent random variables such that
$Z_1\sim \mathrm{gamma}(\alpha,1)$ and $Z_2\sim \mathrm{gamma}(\beta,1)$. Since $Z_1+Z_2\sim \mathrm{gamma}(\alpha+\beta)$ is independent of $\bab$, the relation $X\eqdis Y\bab$ is equivalent to $X(Z_1+Z_2) \eqdis YZ_1$. It follows from Jessen  and Mikosch (2006, Lemma 4.2(a)) that  the survivor function of $Y Z_1$ is regularly varying with index $-\gamma$, and Lemma 17 in Hashorva et al. (2007) implies the same is true for $\o H(x)$. We emphasize that this proof is valid for $\gamma\ge0$.  \QED

Note that Theorem 12.3.2 in Berman (1992) follows from the above direct proof since
$$ \E{(1-\bab)^{\gamma/2}}=\E{ B_{\beta,\alpha}^{\gamma/2}}={B(\alpha,\beta+\gamma/2) \over B(\alpha,\beta)}.$$

The situation where $B$ is allowed to be unbounded can be handled by writing
\BQN\label{eq:decomp}
\pk{YB>x}= \pk{YB>x; Y>x}+\pk{YB>x; Y\le x}.
\EQN
Exactly as in the last proof, the first term on the right is asymptotically proportional to $\o H(x)\pk{B> W^{-1}}$, but now the probability term evaluates as
$$\pk{W> B^{-1}}=\pk{B\ge 1}+\E{B^\gamma; B\le 1}.$$
The second term on the right-hand side of \eqref{eq:decomp} equals
\begin{eqnarray*}
\pk{x/B \le Y\le x, B>1} &=& \int_1^\infty \left(\o H(x/z)-\o H(x)\right) d \pk{B\le z}\\
&=& (1+o(1)) \o H(x)\left[\E{B^\gamma; B>1}-\pk{B>1}\right],
\end{eqnarray*}
provided the limit here can be taken inside the integral. This is permissible if $\E{B^{\gamma+\epsilon}}<\infty$ for some $\epsilon>0$.

The converse tail equivalence statement is open in general, but see Hashorva et al. (2007) for the case where $B$ has a gamma distribution.

The following result is the analogue of \netheo{theo:new:G}   for $H\in MDA(\Psi_\gamma)$ and it generalizes the direct assertion of \netheo{eq:theo:BM1:3}.
{\BT \label{theo:new:W} Let $H(0)=0$,  $r_H=1$ and  $H\in MDA(\Psi_\gamma)$. (a) If \eqref{eq:betaTail} holds, then
$$\hat I(x)=\int_{1-x}^1 \phi(x/y) dH(y)= (1+o(1)) C{\Gamma(\beta+1)\Gamma(\gamma+1) \over \Gamma(\beta+\gamma+1)}
x^\beta \o H(1-x),\quad(x\downarrow 0).$$
(b) If \eqref{eq:betaTail:b} holds, then
$$\hat J(x)=\int_{1-x}^1 y^{-1}g(x/y)dH(y) = (1+o(1)) c {\Gamma(\beta)\Gamma(\gamma+1) \over \Gamma(\beta+\gamma)}
x^{\beta-1} \o H(1-x),\quad(x\downarrow 0).$$
\ET}

\prooftheo{theo:new:W} For (a) simply observe that if $1-x<y<1$, then
$$\phi\left(1-{y-x \over y}\right)= (1+o(1)) C(y-x)^{\beta} y^{-\beta},\quad(x \downarrow 0).$$
Hence $\hat I(x)$ is asymptotically equal to $C \Gamma(\beta+1) (\OPJ_{\beta+1,p_{-\beta}}H)(1-x)$, and the assertion follows from \netheo{lem:weib:1}. Similarly, $\hat J(x)$ is asymptotically equal to $c \Gamma(\beta) (\OPJ_{\beta,p_{-\beta}}H)(1-x)$.\QED

\prooftheo{eq:theo:BM1:3} If $H\in MDA(\Psi_\gamma)$, then  \eqref{res:BERM3} and \eqref{res:BERM3:d} follow from \netheo{theo:new:W}.

\prooflem{power} It follows from \eqref{eq:uv:B} that
$$(y+t/w(y))^p = (1+o(1)) py^{p-1}/w(y) =(1+o(1))\left(y^p+(t/w_p(y^p))\right),$$
and hence that the necessary and sufficient condition \eqref{eq:gumbel} applied to $F$ is equivalent to
$$\lim_{y\to\infty}\pk{X^p>y^p+t/w_p(y^p) | X^p>y^p} = e^{-t}.$$
Setting $x=y^p$ shows this is equivalent to $F_p\in MDA(\Lambda,w_p)$. \QED

\prooftheo{eq:ellip:A}  Let $G_2$ and $H_2$ denote the distribution functions of $U^2$ and $R^2$, respectively. It follows from \eqref{eee},  \nelem{power} and \netheo{eq:theo:BM1:0} that
$$H\in MDA(\Lambda,w) \ {\rm iff} \ H_2\in MDA(\Lambda,w_2) \ {\rm iff} \
G_2\in MDA(\Lambda,w_2) \ {\rm iff} \  G\in MDA(\Lambda,w).$$
This, together with Theorem 12.3.3 in Berman (1992)  implies that  \eqref{ajm:2} holds if $G\in MDA(\Lambda,w)$, i.e. (a) is valid.

By the same reasoning, (b) follows if we prove it  assuming $H\in MDA(\Lambda,w)$ and $H$ has a density  function $h$. Observing that $Z_x$ has the same distribution as $\sqrt{1-\rho^2} S_2|S_1=x$, we set $\rho=0$ without loss of generality. In this case (following the example of Abdous et al. (2005))  we can use the equivalent representation $(U,V)=(I_1\sqrt B,I_2\sqrt{(1-B^2})$, where $I_1$, $I_2$ and $B$ are independent, the $I_j=\pm 1$ with equal probability, and $B\sim \mathrm{beta}(1/2,1/2)$. The joint density function $f(u,v)$ of $(U,v)$ is radially symmetric and a routine computation yields
$$f(u,v)=\left[2\pi\sqrt{u^2+v^2}\right]^{-1}h\left(\sqrt{u^2+v^2}\right).$$
It is more expedient to work directly in terms of $h_2(z)=\left(2\sqrt z\right)^{-1}h\left(\sqrt z\right)$, whence
$$f(u,v)=\pi^{-1}h_2\left(u^2+v^2\right).$$
Integration with respect to $v$ and using the substitution $y=v^2$ gives the marginal density function of $U$,
$$f_U(u)={1\over \pi}\int_0^\infty h_2(y+u^2)y^{-1/2}dy,$$
and hence the density function of $Z_x$ is
$$f(v|x):=f(x,v)/f_U(x)={h_2\left(x^2+v^2\right) \over \int_0^\infty h_2(y+v^2)y^{-1/2}dy},$$
valid for real $v$ and $x>0$. Note that the distribution of $Z_x$ is symmetric about zero.

Let $t>0$ and replace $v$ with $t/c(x)$ in this density function. Since $c^2(x)=2w_2(x^2)$, the density function of $c(x)Z_x$ is the function of $s=x^2$ given by
\BQN\label{nd}
\zeta(t|x)={h_2(s+t^2/2w_2(s)) \over \sqrt{2w_2(s)}\int_0^\infty h_2(s+y)y^{-1/2}dy}.
\EQN
Divide the numerator and denominator of the right-hand side by $w_2(s) \o H_2(s)$.  Since $h_2(s)=(1+o(1))w_2(s)\o H_2(s)$, it follows from \nelem{power}, and  \eqref{eq:gumbel} and \eqref{eq:self} applied to $H_2$, that the numerator term obtained from \eqref{nd} converges to $\exp(-t^2/2)$, as $x\to r_h$.

Next, making the substitution $z=yw_2(s)$ in the integral at \eqref{nd}, the denominator term obtained from the division operation is
$${\sqrt 2 \over w_2(s)\o H_2(s)}\int_0^\infty h_2(s+z/w_2(s)) z^{-1/2} dz=\sqrt 2 \vek{E}\left\{W_s^{-1/2}\right\},$$
where $W_s$ is as defined in the proof of \netheo{theo:OPJ:G} ($s$ replacing $x$ there). The moment convergence theorem ensures that
$$\lim_{x\to r_H} \vek{E}\left\{W_s^{-1/2}\right\}=\vek{E}\left\{W^{-1/2}\right\} =\Gamma(1/2)=\sqrt\pi.$$
Combining these limits shows that $ \zeta(t|x)$ converges to the standard Gaussian density function, and the assertion follows.
\QED



\bibliographystyle{plain}

\end{document}